\documentclass{article}
\usepackage{amsfonts}

\newtheorem{Theorem}{\bf Theorem}[section]
\newtheorem{Proposition}[Theorem]{\bf Proposition}
\newtheorem{Corollary}[Theorem]{\bf Corollary}
\newtheorem{Definition}[Theorem]{\bf Definition}

\newtheorem{Remark}[Theorem]{\bf Remark}
\input{tcilatex}

\begin{document}

\title{COVARIANT\ COMPLETELY\ POSITIVE\ LINEAR\ MAPS\ BETWEEN\ LOCALLY\ \ $C^{*}$%
-ALGEBRAS}
\author{MARIA\ JOI\c{T}A\thanks{%
 }\thanks{\textit{This research was supported by grant }CNCSIS-code
A1065/2006.\ \ }}
\maketitle

\begin{abstract}
We prove a covariant version of the KSGNS\ (Kasparov, Stinespring,
Gel'fand,Naimark,Segal) construction for completely positive linear maps
between locally $C^{*}$-algebras. As an application of this construction, we
show that a covariant completely positive linear map $\rho $ from a locally $%
C^{*}$-algebra $A$ to another locally $C^{*}$ -algebra $B$ with respect to a
locally $C^{*}$-dynamical system $(G,A,\alpha )$ extends to a completely
positive linear map on the crossed product $A\times _{\alpha ^{{}}}G$.
\end{abstract}

\section{Introduction}

Locally $C^{*}$-algebras are generalizations of $C^{*}$-algebras. Instead of
being given by a single norm, the topology on a locally $C^{*}$-algebra is
defined by a directed family of $C^{*}$-seminorms. Such important concepts
as Hilbert $C^{*}$-modules, adjointable operators, (completely) positive
linear maps, $C^{*}$-dynamical systems can be defined with obvious
modifications in the framework of locally $C^{*}$-algebras. The proofs are
not always straightforward.

It is well-known that a positive linear functional on a $C^{*}$-algebra $A$
induces a representation of this $C^{*}$-algebra on a Hilbert space by the
GNS (Gel'fand, Naimark, Segal) construction (see, for example, [\textbf{1}])
. Stinespring [\textbf{13}] extends this construction for completely
positive linear map from $A$ to $L(H)$, the $C^{*}$-algebra of all bounded
linear operators on a Hilbert space $H.$ On the other hand, Paschke [\textbf{%
8}] (respectively, Kasparov [\textbf{5}]) shows that a completely positive
linear map from $A$ to another $C^{*}$-algebra $B\;$(respectively, from $A$
to the $C^{*}$-algebra of all adjointable operators on the Hilbert $C^{*}$%
-module $H_B$ ) induces a representation of $A$ on a Hilbert $B$ -module. In
[\textbf{2}], the author extends the KSGNS ( Kasparov, Stinespring,
Gel'fand, Naimark, Segal) construction for a strict continuous, completely
positive linear map from a locally $C^{*}$-algebra $A$ to $L_B(E)$, the
locally $C^{*}$-algebra of all adjointable operators on a Hilbert module $E$
over a locally $C^{*}$-algebra $B$. In this paper we propose to prove a
covariant version of this construction. Thus we show that a covariant
completely positive linear map from $A$ to $L_B(E)$ with respect to a
locally $C^{*}$-dynamical system $(G,A,\alpha )$ induces a non-degenerate,
covariant representation of $(G,A,\alpha )$ on a Hilbert $B$ -module which
is unique up to unitary equivalence, Theorem 3.6. Using the analog of the
covariant version of Stinespring construction [\textbf{9}] for bounded
operators on Hilbert $C^{*}$-modules, Kaplan [\textbf{4}] shows that a
discrete covariant completely positive map $\rho $ from a unital $C^{*}$%
-algebra $A$ to another unital $C^{*}$-algebra $B$ extends to a completely
positive map from the crossed product $A\times _{\alpha _{}}G$ to $B$. We
extend this result showing that a non-degenerate, covariant, continuous
completely positive linear map from a locally $C^{*}$-algebra to another
locally $C^{*}$-algebra $B$ extends to a non-degenerate, continuous
completely positive linear map on the crossed product $A\times _{\alpha
_{}}G $, Proposition 3.9.

\section{Preliminaries}

A locally $C^{*}$-algebra is a complete complex Hausdorff topological $*$
-algebra whose topology is determined by a directed family of $C^{*}$%
-seminorms. If $A$ is a locally $C^{*}$-algebra and $S(A)$ is the set of all
continuous $C^{*}$-seminorms on $A,$ then for each $p\in S(A),$ $A_p=A/\ker
p $ is a $C^{*}$-algebra in the norm induced by $p$, and $\{A_p;\pi
_{pq}\}_{p,q\in S(A),p\geq q},$ where $\pi _{pq}$ is the canonical morphism
from $A_p$ onto $A_q$ defined by $\pi _{pq}\left( a+\ker p\right) =a+\ker q$
for all $a\in A$ is an inverse system of $C^{*}$-algebras. Moreover, $A$ can
be identified with $\lim\limits_{\leftarrow p}A_p.$ The canonical map from $%
A $ onto $A_p$ is denoted by $\pi _p.$

An approximate unit of $A$ is an increasing net $\{e_{\lambda
^{}}\}_{\lambda \in \Lambda }$ of positive elements of $A$ such that $%
p(e_{\lambda ^{}})\leq 1$ for all $p\in S(A)$ and for all $\lambda \in
\Lambda ,$ $p(ae_{\lambda ^{}}-a)\rightarrow 0$ and $p(e_{\lambda
^{}}a-a)\rightarrow 0$ for all $p\in S(A)$ and for all $a\in A.$ Any locally 
$C^{*}$-algebra has an approximate unit [\textbf{11}, Proposition 3.11].

A morphism of locally $C^{*}$-algebras is a continuous $*$ -morphism from a
locally $C^{*}$-algebra $A$ to another locally $C^{*}$-algebra $B$. An
isomorphism of locally $C^{*}$-algebras from $A$ to $B$ is a bijective map $%
\Phi :$ $A$ $\rightarrow $ $B$ such that $\Phi $ and $\Phi ^{-1}$ are
morphisms of locally $C^{*}$-algebras.

Let $M_n(A)$ denote the $*$ -algebra of all $n\times n$ matrices over $A$
with the algebraic operations and the topology obtained by regarding it as a
direct sum of $n^2$ copies of $A.$ Then $\{M_n(A_p);\pi
_{pq}^{(n)}\}_{p,q\in S(A),p\geq q},$ where $\pi _{pq}^{(n)}([\pi
_p(a_{ij})]_{i,j=1}^n)=$ $[\pi _q(a_{ij})]_{i,j=1}^n,$ is an inverse system
of $C^{*}$-algebras and $M_n(A)$ can be identified with $\lim\limits_{%
\leftarrow p}M_n(A_p).$

A linear map $\rho :A\rightarrow B$ between two locally $C^{*}$-algebras is
completely positive if the linear maps $\rho ^{(n)}:M_n(A)\rightarrow M_n(B)$
defined by 
\[
\rho ^{(n)}([a_{ij}]_{i,j=1}^n)=[\rho (a_{ij})]_{i,j=1}^n 
\]
$n=1,2,...,n,...,$ are all positive.

\begin{Definition}
\emph{\ A }pre-Hilbert$\ A$-module\emph{\ is a complex vector space}$\ E$%
\emph{\ which is also a right }$A$\emph{-module, compatible with the complex
algebra structure, equipped with an }$A$\emph{-valued inner product }$%
\left\langle \cdot ,\cdot \right\rangle :E\times E\rightarrow A\;$\emph{%
which is }$\Bbb{C}$\emph{\ -and }$A$\emph{-linear in its second variable and
satisfies the following relations:}

\begin{enumerate}
\item 
\begin{enumerate}
\item 
\begin{enumerate}
\item  $\;\left\langle \xi ,\eta \right\rangle ^{*}=\left\langle \eta ,\xi
\right\rangle \;\;$\emph{for every }$\xi ,\eta \in E;$

\item  $\;\left\langle \xi ,\xi \right\rangle \geq 0\;\;$\emph{for every }$%
\xi \in E;$

\item  $\left\langle \xi ,\xi \right\rangle =0\;$\emph{\ if and only if }$%
\xi =0.$
\end{enumerate}
\end{enumerate}
\end{enumerate}

\emph{We say that }$E\;$\emph{is a Hilbert }$A$\emph{-module if }$E\;$\emph{%
is complete with respect to the topology determined by the family of
seminorms }$\{\left\| \cdot \right\| _{p}\}_{p\in S(A)}\;$\emph{where }$%
\left\| \xi \right\| _{p}=\sqrt{p\left( \left\langle \xi ,\xi \right\rangle
\right) },\xi \in E$\emph{\ [11, Definition 4.1].}
\end{Definition}

\smallskip

Let $E\;$be a Hilbert $A$-module.\ For $p\in S(A),\;\mathcal{E}_p=\{\xi \in
E;p(\left\langle \xi ,\xi \right\rangle )=0\}\;$is a closed submodule of $%
E\; $and $E_p=E/\mathcal{E}_p\;$is a Hilbert $A_p$-module with $(\xi +%
\mathcal{E}_p)\pi _p(a)=\xi a+\mathcal{E}_p\;$and $\left\langle \xi +%
\mathcal{E}_p,\eta +\mathcal{E}_p\right\rangle =\pi _p(\left\langle \xi
,\eta \right\rangle ).$\ The canonical map from $E\;$onto $E_p$ is denoted
by $\sigma _p.$ For $p,q\in S(A),\;p\geq q\;$there is a canonical morphism
of vector spaces $\sigma _{pq}\;$from $E_p\;$onto $E_q\;$such that $\sigma
_{pq}(\sigma _p(\xi ))=\sigma _q(\xi ),\;\xi \in E.\;$Then $\{E_p;A_p;\sigma
_{pq}\}_{p,q\in S(A),p\geq q}$ is an inverse system of Hilbert $C^{*}$%
-modules in the following sense: $\sigma _{pq}(\xi _pa_p)=\sigma _{pq}(\xi
_p)\pi _{pq}(a_p),\xi _p\in E_p,a_p\in A_p;$ $\left\langle \sigma _{pq}(\xi
_p),\sigma _{pq}(\eta _p)\right\rangle =\pi _{pq}(\left\langle \xi _p,\eta
_p\right\rangle ),\xi _p,\eta _p\in E_p;$ $\sigma _{pp}(\xi _p)=\xi _p,\;\xi
_p\in E_p\;$and $\sigma _{qr}\circ \sigma _{pq}=\sigma _{pr}\;$if $p\geq
q\geq r,\ $and $\lim\limits_{\leftarrow p}E_p$ is a Hilbert $A$-module which
can be identified with\ $E$ [\textbf{11}, Proposition 4.4].

Let $E\;$and $F$\ be Hilbert $A$-modules. We say that an $A$-module morphism 
$T:E\rightarrow F\;$is adjointable if there is an $A$-module morphism $%
T^{*}:F\rightarrow E$\ such that $\left\langle T\xi ,\eta \right\rangle
=\left\langle \xi ,T^{*}\eta \right\rangle \;$for every $\xi \in E$ and $%
\eta \in F.$ Any adjointable $A$-module morphism is continuous. The set $%
L_{A}(E,F)$ of all adjointable $A$-module morphisms from $E$ into $F$
becomes a locally convex space with topology defined by the family of
seminorms $\{\widetilde{p}\}_{p\in S(A)},$ where $\widetilde{p}(T)=\left\|
(\pi _{p})_{*}(T)\right\| _{L_{A_{p}}(E_{p},F_{p})},$ $T\in L_{A}(E,F)\;$and 
$(\pi _{p})_{*}(T)(\xi +\mathcal{E}_{p})=T\xi +\mathcal{F}_{p},$ $\xi \in E.$
Moreover, $\{L_{A_{p}}(E_{p},F_{p});\;(\pi _{pq})_{*}$ $\}_{p,q\in
S(A),p\geq q},$ where $(\pi _{pq})_{*}:L_{A_{p}}(E_{p},F_{p})\rightarrow
L_{A_{q}}(E_{q},F_{q}),$ $(\pi _{pq})_{*}(T_{p})(\sigma _{q}(\xi ))=\chi
_{pq}(T_{p}(\sigma _{p}(\xi ))),$ and $\chi _{pq},\;p,q\in S(A),$ $p\geq q\;$%
are the connecting maps of the inverse system $\{F_{p}\}_{p\in S(A)},$ is an
inverse system of Banach spaces, and $\lim\limits_{\leftarrow
p}L_{A_{p}}(E_{p},F_{p})$ can be identified with $L_{A}(E,F)$ [\textbf{11},
Proposition 4.7]. Thus topologized, $L_{A}(E,E)$ becomes a locally $C^{*}$%
-algebra, and we write $L_{A}(E)\;$for $L_{A}(E,E).$

The strict topology on $L_A(E)$ is defined by the family of seminorms $%
\{\left\| \cdot \right\| _{p,\xi }$ $\}_{(p,\xi )\in S(A)\times E},$ where $%
\left\| T\right\| _{p,\xi }=\left\| T\xi \right\| _p+\left\| T^{*}\xi
\right\| _p,$ $T\in L_A(E).$

Two Hilbert $A$ -modules $E$ and $F$ are unitarily equivalent if there is a
unitary element in $L_A(E,F).$

A non-degenerate representation of a locally $C^{*}$-algebra $A$ on a
Hilbert module $E$ over a locally $C^{*}$-algebra $B$ is a morphism of
locally $C^{*}$-algebras $\Phi $ from $A$ to $L_B(E)\;$such that $\Phi
\left( A\right) E$ is dense in $E.$

A continuous completely positive linear map $\rho $ from $A$ to $L_B(E)$ is
non-degenerate if the net $\{\rho \left( e_{\lambda _{}}\right) \}_{\lambda
\in \Lambda }$ converges strictly to the identity map on $E$ for some
approximate unit $\{e_{\lambda _{}}\}_{\lambda \in \Lambda }$ for $A.$

Let $G$ be a locally compact group and let $A$ be a locally $C^{*}$-algebra.
An action of $G$ on $A$ is a morphism $\alpha $ from $G$ to Aut$\left(
A\right) $, the set of all isomorphisms of locally $C^{*}$-algebras from $A$
to $A$. The action $\alpha $ is continuous if the function $\left(
t,a\right) \rightarrow \alpha _t(a)$ from $G\times A$ to $A$ is jointly
continuous. An action $\alpha $ is called an inverse limit action if we can
write $A$ as inverse limit $\lim\limits_{\leftarrow \delta }A_{\delta _{}}$
of $C^{*}$-algebras in such a way that there are actions $\alpha ^{\left(
\delta \right) }$ of $G$ on $A_{\delta _{}}$ such that $\alpha
_t=\lim\limits_{\leftarrow \delta }\alpha _t^{\left( \delta \right) }$ for
all $t$ in $G\ $[\textbf{12}, Definition 5.1]. An action $\alpha $ of $G$ on 
$A$ is a continuous inverse limit action if there is a cofinal subset $%
S_G(A,\alpha )$ of $G$-invariant continuous $C^{*}$-seminorms on $A$ ( a
continuous $C^{*}$-seminorm $p$ on $A$ is $G$ -invariant if $p(\alpha
_t(a))=p(a)$ for all $a$ in $A$ and for all $t$ in $G$). So if $\alpha $ is
a continuous inverse limit action of $G$ on $A$ we can suppose that $%
S(A)=S_G(A,\alpha ).$

A locally $C^{*}$-dynamical system is a triple $(G,A,\alpha ),$ where $G$ is
a locally compact group, $A$ is a locally $C^{*}$-algebra and $\alpha $ is a
continuous action of $G$ on $A.$

Let $\alpha $ be a continuous inverse limit action of $G$ on $A.$ The set $%
C_c(G,A)$ of all continuous functions from $G$ to $A$ with compact support
becomes a $*$ -algebra with convolution of two functions 
\[
\left( f\times h\right) \left( s\right) =\int\limits_Gf(t)\alpha _t\left(
h(t^{-1}s)\right) dt 
\]
as product and involution defined by 
\[
f^{\sharp }(t)=\Delta (t)^{-1}\alpha _t\left( f(t^{-1})^{*}\right) 
\]
where $\Delta $ is the modular function on $G.$ The Hausdorff completion of $%
C_c(G,A)$ with respect to the topology defined by the family of
submultiplicative $*$ -seminorms $\{N_p\}_{p\in S(A)},$ where 
\[
N_p(f)=\int\limits_Gp(f(s))ds 
\]
is denoted by $L^1(G,A,\alpha )$ and the enveloping locally $C^{*}$-algebra $%
A\times _{\alpha ^{}}G$ of $L^1(G,A,\alpha )$ is called the crossed product
of $A$ by $\alpha $ [\textbf{3}, Definition 3.14]. Moreover, the $C^{*}$%
-algebras $(A\times _{\alpha ^{}}G)_p$ and $A_p\times _{\alpha ^{(p)}}G$ can
be identified for each $p\in S(A)$ and so $A\times _{\alpha ^{}}G$ can be
identified with $\lim\limits_{\leftarrow p}$ $A_p\times _{\alpha ^{(p)}}G$ [%
\textbf{3}, Remark 3.15].

\section{Covariant representations associated with a covariant completely
positive linear map}

Let $B$ be a locally $C^{*}$-algebra, let $E$ be a Hilbert $B$ -module and
let $G$ be a locally compact group.

\begin{Definition}
\emph{A }unitary representation \emph{of }$G$\emph{\ on }$E$\emph{\ is a map 
}$u$\emph{\ from }$G$\emph{\ to }$L_{B}(E)$\emph{\ such that}

\begin{enumerate}
\item 
\begin{enumerate}
\item  $u_{g}$\emph{\ is a unitary element in }$L_{B}(E)$\emph{\ for all }$%
g\in G;$

\item  $u_{gt}=u_{g}u_{t}$\emph{\ for all }$g,t\in G;$

\item  \emph{the map }$g\mapsto u_{g}\xi $\emph{\ from }$G$\emph{\ to }$E$%
\emph{\ is continuous for all }$\xi \in E.$
\end{enumerate}
\end{enumerate}
\end{Definition}

\smallskip

\begin{Remark}
If $u$\ is a unitary representation of $G$\ on $E,$\ then for each $q\in
S(B),$\ $g\mapsto (\pi _{q})_{*}\circ u$\ is a unitary representation of $G$%
\ on $E_{q}.$\ Moreover, $u_{g}=\lim\limits_{\leftarrow q}u_{g}^{(q)},$\
where $u_{g}^{(q)}=(\pi _{q})_{*}(u_{g}),$\ for all $g\in G.$
\end{Remark}

\smallskip

\begin{Definition}
A non-degenerate, covariant representation\emph{\ of a locally }$C^{*}$\emph{%
-dynamical system }$(G,A,\alpha )$\emph{\ on a Hilbert }$B$\emph{\ -module }$%
E$\emph{\ is a triple }$(\Phi ,v,E)$\emph{, where }$\Phi $\emph{\ is a
non-degenerate representation of }$A$\emph{\ on }$E$\emph{, }$v$\emph{\ is a
unitary representation of }$G$\emph{\ on }$E$\emph{\ and } 
\[
\Phi (\alpha _{g}(a))=v_{g}\Phi (a)v_{g}^{*} 
\]
\emph{for all }$g\in G$\emph{\ and }$a\in A.$
\end{Definition}

\emph{\smallskip }

\begin{Proposition}
Let $(G,A,\alpha )$ be a locally $C^{*}$-dynamical system such that $\alpha $
is an inverse limit action, let $(\Phi ,v,E)$ be a non-degenerate covariant
representation of $(G,A,\alpha )$ on a Hilbert $B$ -module $E$. Then there
is a unique non-degenerate representation $\Phi \times v$ of the crossed
product $A\times _{\alpha ^{{}}}G$ on $E$ such that 
\[
(\Phi \times v)(f)=\int\limits_{G}\Phi (f(g))v_{g}dg 
\]
for all $f\in C_{c}(G,A).$
\end{Proposition}

PROOF. We partition the proof into two steps.

\textit{Step1. }We suppose that $B$ is a $C^{*}$-algebra.

Since $\Phi $ is a representation of $A$ on $E$, there is $p\in S(A)$ such
that $\left\| \Phi \left( a\right) \right\| _{L_B(E)}\leq p(a)$ for all $%
a\in A$. From this fact, we deduce that there is a morphism of $C^{*}$%
-algebras $\Phi _p$ from $A_p$ to $L_B(E)$ such that $\Phi _p\circ \pi
_p=\Phi .$ Therefore $\Phi _p$ is a representation of $A_p$ on $E,$ and
moreover, it is non-degenerate, since $\Phi $ is non-degenerate and $\pi _p$
is surjective. It is not difficult to check that $(\Phi _p,v,E)$ is a
non-degenerate covariant representation of $(G,A_p,\alpha ^{(p)}).$ Then
there is a unique non-degenerate representation $\Phi _p\times v$ of $%
A_p\times _{\alpha ^{(p)}}G$ on $E$ such that 
\[
(\Phi _p\times v)(f)=\int\limits_G\Phi _p(f(g))v_gdg 
\]
for all $f$ $\in C_c(G,A_p)$ (see, for example, Proposition 7.6.4, [\textbf{%
10}]). Therefore $\Phi \times v=(\Phi _p\times v)\circ \widetilde{\pi }_p,$
where $\widetilde{\pi }_p$ is the canonical map from $A\times _{\alpha ^{}}G$
onto $A_p\times _{\alpha ^{(p)}}G$ is a non-degenerate representation of $%
A\times _{\alpha ^{}}G$ on $E$ such that 
\[
(\Phi \times v)(f)=\int\limits_G\Phi _p(\widetilde{\pi }_p(f)(g))v_gdg=\int%
\limits_G\Phi (f(g))v_gdg 
\]
for all $f\in C_c(G,A),$ and since $C_c(G,A)$ is dense in $A\times _{\alpha
^{}}G,$ $\Phi \times v$ is unique with the above property.

\textit{Step2.} The general case.

For each $q\in S(B),$ $(\pi _q)_{*}\circ \Phi $ is a non-degenerate
representation of $A$ on $E_q,$ and $((\pi _q)_{*}\circ \Phi ,v^{(q)},E_q)$
is a non-degenerate covariant representation of $(G,A,\alpha )$ on $E_q.$ By
Step 1 there is a unique non-degenerate representation $\left( (\pi
_q)_{*}\circ \Phi \right) \times v^{(q)}$ of $A\times _{\alpha ^{}}G$ on $%
E_q $ such that 
\[
(\left( (\pi _q)_{*}\circ \Phi \right) \times v^{(q)})(f)=\int\limits_G(\pi
_q)_{*}(\Phi (f(g)))v_g^{(q)}dg 
\]
for all $f\in C_c(G,A).$ By Lemma 3.7 in [\textbf{3}], we have 
\begin{eqnarray*}
(\pi _{qr})_{*}((\left( (\pi _q)_{*}\circ \Phi \right) \times v^{(q)})(f))
&=&\int\limits_G(\pi _{qr})_{*}((\pi _q)_{*}(\Phi (f(g)))v_g^{(q)})dg \\
&=&\int\limits_G(\pi _r)_{*}(\Phi (f(g)))v_g^{(r)})dg \\
&=&(\left( (\pi _r)_{*}\circ \Phi \right) \times v^{(r)})(f)
\end{eqnarray*}
for all $f\in C_c(G,A)$ and for all $q,r\in S(B)$ with $q\geq r.$ Therefore $%
(\pi _{qr})_{*}\circ (\left( (\pi _q)_{*}\circ \Phi \right) \times
v^{(q)})=\left( (\pi _r)_{*}\circ \Phi \right) \times v^{(r)}$ for all $%
q,r\in S(B)$ with $q\geq r.$ This implies that there is a continuous $*$
-morphism $\Phi \times v$ from $A\times _{\alpha ^{}}G$ to $L_B(E)$ such
that $(\pi _q)_{*}\circ (\Phi \times v)=$ $\left( (\pi _q)_{*}\circ \Phi
\right) \times v^{(q)}$ for all $q\in S(B).$ Using Lemma III 3.1 in [\textbf{%
7}], it is not hard to check that $(\Phi \times v)(A\times _{\alpha ^{}}G)E$
is dense in $E.$ Therefore $\Phi \times v$ is a non-degenerate
representation of $A\times _{\alpha ^{}}G$ on $E$, and by Lemma 3.7 in [%
\textbf{3}], 
\[
(\Phi \times v)(f)=\int\limits_G\Phi (f(g))v_gdg 
\]
for all $f\in C_c(G,A).$ Moreover, since $C_c(G,A)$ is dense in $A\times
_{\alpha ^{}}G,$ $\Phi \times v$ is unique with the above property. q.e.d.

\smallskip

\begin{Definition}
\emph{Let }$(G,A,\alpha )$\emph{\ be a locally }$C^{*}$\emph{-dynamical
system and let }$u$\emph{\ be a unitary representation of }$G$\emph{\ on a
Hilbert }$B$\emph{\ -module }$E$\emph{. We say that a completely positive
linear map }$\rho $\emph{\ from }$A$\emph{\ to }$L_{B}(E)$\emph{\ is }$u$%
\emph{\ }-covariant\emph{\ with respect to the locally }$C^{*}$\emph{%
-dynamical system }$(G,A,\alpha )$\emph{\ if } 
\[
\rho (\alpha _{g}(a))=u_{g}\rho (a)u_{g}^{*} 
\]
\emph{for all }$a\in A$\emph{\ and for all }$g\in G.$
\end{Definition}

\smallskip

Recall that if $\rho $ is a completely positive linear map from a $C^{*}$%
-algebra $A$ to $L_B(E),$ the $C^{*}$-algebra of all adjointable operators
on a Hilbert module $E$ over a $C^{*}$-algebra $B,$ the quotient vector
space $\left( A\otimes _{\text{alg}}E\right) /\mathcal{N}_{\rho _{}},$ where 
$\mathcal{N}_{\rho _{}}=\{\sum_{i=1}^na_i\otimes \xi _i;$ $%
\sum_{i,j=1}^n\left\langle \xi _i,\rho \left( a_i^{*}a_j\right) \xi
_j\right\rangle =0\}$ becomes a pre-Hilbert $B$ -module with the action of $%
B $ on $A\otimes _{\text{alg}}E$ defined by $\left( a\otimes \xi +\mathcal{N}%
_{\rho _{}}\right) b=a\otimes \xi b+\mathcal{N}_{\rho _{}}$ and the
inner-product defined by 
\[
\left\langle \sum\limits_{i=1}^na_i\otimes \xi _i+\mathcal{N}_{\rho
_{}},\sum\limits_{j=1}^mb_j\otimes \eta _j+\mathcal{N}_{\rho
_{}}\right\rangle _{\rho
_{}}=\sum\limits_{i=1}^n\sum\limits_{j=1}^m\left\langle \xi _i,\rho \left(
a_i^{*}b_j\right) \eta _j\right\rangle . 
\]

The following theorem is a covariant version of Theorem 4.6 in [\textbf{2}].

\begin{Theorem}
Let $(G,A,\alpha )$ be a locally $C^{*}$-dynamical system, let $u$ be a
unitary representation of $G$ on a Hilbert module $E$ over a locally $C^{*}$%
-algebra $B$, and let $\rho $ be a $u$ -covariant, non-degenerate,
continuous completely positive linear map from $A$ to $L_{B}(E).$

\begin{enumerate}
\item  Then there is a covariant representation $(\Phi _{\rho _{{}}},v^{\rho
_{{}}},E_{\rho _{{}}})$ of $(G,A,\alpha )$ and an element $V_{\rho _{{}}}$
in $L_{B}(E,E_{\rho _{{}}})$ such that

\begin{enumerate}
\item  $\rho (a)=$ $V_{\rho _{{}}}^{*}\Phi _{\rho _{{}}}(a)V_{\rho _{{}}}$
for all $a\in A;$

\item  $\{\Phi _{\rho _{{}}}(a)V_{\rho _{{}}}\xi ;a\in A,\xi \in E\}$ spans
a dense submodule of $E_{\rho _{{}}};$

\item  $v_{g}^{\rho _{{}}}V_{\rho _{{}}}=V_{\rho _{{}}}u_{g}$ for all $g\in
G.$
\end{enumerate}

\item  If $F$ is a Hilbert $B$ -module, $(\Phi ,v,F)$ is a covariant
representation of $(G,A,\alpha )$ and $W$ is an element in $L_{B}(E,F)$ such
that

\begin{enumerate}
\item  $\rho (a)=W^{*}\Phi (a)W$ for all $a\in A;$

\item  $\{\Phi (a)W\xi ;a\in A,\xi \in F\}$ spans a dense submodule of $F;$

\item  $v_{g}W=Wu_{g}$ for all $g\in G,$
\end{enumerate}
\end{enumerate}

then there is a unitary operator $U$ in $L_{B}(E_{\rho _{{}}},F)$ such that

\begin{enumerate}
\item 
\begin{enumerate}
\item 
\begin{enumerate}
\item  $\Phi (a)U=U\Phi _{\rho _{{}}}(a)$ for all $a\in A;$

\item  $v_{g}U=Uv_{g}^{\rho _{{}}}$ for all $g\in G;$

\item  $W=UV_{\rho _{{}}}.$
\end{enumerate}
\end{enumerate}
\end{enumerate}
\end{Theorem}

PROOF. We partition the proof into two steps.

\textit{Step 1. }Suppose that $B$ is a $C^{*}$-algebra.

$1.\;$Let $\{e_{\lambda ^{}}\}_{\lambda \in \Lambda }$ be an approximate
unit of $A$ such that the net $\{\rho (e_{\lambda ^{}})\}_{\lambda \in
\Lambda }$ is strictly convergent to the identity operator on $E,$ and let $%
(\Phi _{\rho _{}};V_{\rho _{}};E_{\rho _{}})$ be the KSGNS\ construction
associated with $\rho $. Since $\rho $ is continuous there is $p\in S(A)$
and a completely positive linear map $\rho _p$ from $A_p$ to $L_B(E)$ such
that $\rho =\rho _p\circ \pi _p$ ( see, for example, the proof of
Proposition 3.5 in [\textbf{2}]). By the proof of Theorem 4.6 in [\textbf{2}%
] we can suppose that $E_{\rho _{}}$ is the completion of the pre-Hilbert
space $\left( A_p\otimes _{\text{alg}}E\right) /\mathcal{N}_{\rho _p},$ $%
V_{\rho _{}}\xi =\lim\limits_{\lambda _{}}\left( \pi _p\left( e_{\lambda
_{}}\right) \otimes \xi +\mathcal{N}_{\rho _p}\right) $ and $\Phi _{\rho
_{}}\left( a\right) \left( \pi _p\left( b\right) \otimes \xi +\mathcal{N}%
_{\rho _p}\right) =\pi _p\left( ab\right) \otimes \xi +\mathcal{N}_{\rho _p}$
for all $a,b\in A$ and for all $\xi \in E.$

Let $g\in G.$ From

$\left\langle \sum\limits_{i=1}^n\pi _p\left( a_i\right) \otimes \xi _i+%
\mathcal{N}_{\rho _p},\sum\limits_{j=1}^mb_j\otimes \eta _j+\mathcal{N}%
_{\rho _p}\right\rangle _{\rho
_p}=\sum\limits_{i=1}^n\sum\limits_{j=1}^m\left\langle \xi _i,\rho \left(
a_i^{*}b_j\right) \eta _j\right\rangle $

$\,\,=\sum\limits_{i=1}^n\sum\limits_{j=1}^m\left\langle u_g\left( \xi
_i\right) ,u_g\rho \left( a_i^{*}b_j\right) u_{g^{-1}}u_g\left( \eta
_j\right) \right\rangle $

$\,\,=\sum\limits_{i=1}^n\sum\limits_{j=1}^m\left\langle u_g\left( \xi
_i\right) ,\rho \left( \alpha _g\left( a_i^{*}b_j\right) \right) u_g\left(
\eta _j\right) \right\rangle $

$\,\,=\left\langle \sum\limits_{i=1}^n\pi _p\left( \alpha _g\left(
a_i\right) \right) \otimes u_g\left( \xi _i\right) +\mathcal{N}_{\rho
_p},\sum\limits_{j=1}^m\pi _p\left( \alpha _g\left( b_j\right) \right)
\otimes u_g\left( \eta _j\right) +\mathcal{N}_{\rho _p}\right\rangle _{\rho
_p}$

for all $\xi _1,...,\xi _n,$ $\eta _1,....,\eta _m\in E,$ for all $%
a_1,...,a_n,$ $b_1,....,b_m\in A$ and for all $g\in G,$ we deduce that,
there is a unitary operator $v_g^{\rho _{}}$ in $L_B\left( E_{\rho
_{}}\right) $ such that 
\[
v_g^{\rho _{}}(\pi _p(a)\otimes \xi +\mathcal{N}_{\rho _p})=\pi _p(\alpha
_g(a))\otimes u_g\xi +\mathcal{N}_{\rho _p} 
\]
for all $a\in A$ and for all $\xi \in E.$ It is not difficult to check that
the map $g\mapsto v_g^{\rho _{}}$ from $G$ to $L_B\left( E_{\rho _{}}\right) 
$ is a unitary representation of $G$ on $E_{\rho _{}}.$

To show that $(\Phi _{\rho _{}},v^{\rho _{}},E_{\rho _{}})$ is a covariant
representation of $(G,A,\alpha )$ it remains to prove that $\Phi _{\rho
_{}}(\alpha _g(a))=v_g^{\rho _{}}\Phi _{\rho _{}}(a)v_{g^{-1}}^{\rho _{}}$
for all $g\in G$ and $a\in A.$ Let $g\in G$ and $a\in A.$ We have

\begin{eqnarray*}
(v_g^{\rho _{}}\Phi _{\rho _{}}(a)v_{g^{-1}}^{\rho _{}})(\pi _p(b)\otimes
\xi +\mathcal{N}_{\rho _p}) &=&(v_g^{\rho _{}}\Phi _{\rho _{}}(a))(\pi
_p(\alpha _{g^{-1}}(b))\otimes u_{g^{-1}}\xi +\mathcal{N}_{\rho _p}) \\
&=&v_g^{\rho _{}}(\pi _p(a\alpha _{g^{-1}}(b))\otimes u_{g^{-1}}\xi +%
\mathcal{N}_{\rho _p}) \\
&=&\pi _p(\alpha _g(a)b)\otimes \xi +\mathcal{N}_{\rho _p} \\
&=&(\Phi _{\rho _{}}(\alpha _g(a)))(\pi _p(b)\otimes \xi +\mathcal{N}_{\rho
_p})
\end{eqnarray*}
for all $b\in A$ and for all $\xi \in E.$ Hence $\Phi _{\rho _{}}(\alpha
_g(a))=v_g^{\rho _{}}\Phi _{\rho _{}}(a)v_{g^{-1}}^{\rho _{}}.$

By Theorem 4.6 (1), [\textbf{2}] the conditions $(a)$ and $(b)$ are
verified. To show that the condition $(c)$ is verified, let $\xi \in E$ and $%
g\in G$. Then we have

$\left\| v_g^{\rho _{}}V_{\rho _{}}\xi -V_{\rho _{}}u_g\xi \right\|
^2=\lim\limits_{\lambda \in \Lambda }\left\| v_g^{\rho _{}}\left( \pi
_p\left( e_{\lambda _{}}\right) \otimes \xi +\mathcal{N}_{\rho _p}\right)
-V_{\rho _{}}u_g\xi \right\| ^2$

$=\lim\limits_{\lambda \in \Lambda }\left\| \left\langle \xi ,\rho
(e_{\lambda ^{}}^2)\xi \right\rangle +\left\langle \xi ,\xi \right\rangle
-\left\langle \rho (\alpha _g(e_{\lambda ^{}}))u_g\xi ,u_g\xi \right\rangle
-\left\langle u_g\xi ,\rho (\alpha _g(e_{\lambda ^{}}))u_g\xi \right\rangle
\right\| $

$\leq \lim\limits_{\lambda \in \Lambda }\left\| \left\langle \xi ,\rho
(e_{\lambda ^{}})\xi \right\rangle +\left\langle \xi ,\xi \right\rangle
-\left\langle \rho (e_{\lambda ^{}})\xi ,\xi \right\rangle -\left\langle \xi
,\rho (e_{\lambda ^{}})\xi \right\rangle \right\| $

$=\lim\limits_{\lambda \in \Lambda }\left\| \left\langle \xi -\rho
(e_{\lambda ^{}})\xi ,\xi \right\rangle \right\| =0.$

Therefore the condition $(c)$ is also verified.

$2.$ By Theorem 4.6 (2), [\textbf{2}] there is a unitary operator $U$ in $%
L_B(E_{\rho _{}},F)$ defined by $U(\Phi _{\rho _{}}(a)V_{\rho _{}}\xi )=\Phi
(a)W\xi $ such that $\Phi (a)U=U\Phi _{\rho _{}}(a)$ for all $a\in A,$ and $%
W=UV_{\rho _{}}.$

Let $g\in G.$ From 
\begin{eqnarray*}
(v_gU)(\Phi _{\rho _{}}(a)V_{\rho _{}}\xi ) &=&v_g(\Phi (a)W\xi )=\Phi
(\alpha _g(a))v_gW\xi \\
&=&\Phi (\alpha _g(a))Wu_g\xi =U(\Phi _{\rho _{}}(\alpha _g(a))V_{\rho
_{}}u_g\xi ) \\
&=&U(\Phi _{\rho _{}}(\alpha _g(a))v_g^{\rho _{}}V_{\rho _{}}\xi
)=(Uv_g^{\rho _{}})(\Phi _{\rho _{}}(a)V_{\rho _{}}\xi ).
\end{eqnarray*}
for all $a\in A$ and for all $\xi \in E,$ we conclude that $v_gU=Uv_g^{\rho
_{}}$ and thus the assertion $2.$ is proved.

\textit{Step 2.} The general case.

Let $q\in S(B).$ Then $\rho _q=(\pi _q)_{*}\circ \rho $ is a $u^{(q)}$
-covariant, non-degenerate, continuous completely positive linear map from $%
A $ to $L_{B_q}(E_q)$, $((\pi _q)_{*}\circ \Phi ,v^{(q)},F_q)$ is a
covariant representation of $(G,A,\alpha )$ and $(\pi _q)_{*}(W)$ is an
element in $L_{B_q}(E_q,F_q)$ such that the conditions $(a),$ $(b)$ and $(c)$
from $2.$ are verified. By Step 1, there is a covariant representation $%
(\Phi _{\rho _q},v^{\rho _q},E_{\rho _q})$ of $(G,A,\alpha )$ and an element 
$V_{\rho _q}$ in $L_{B_q}(E_q,E_{\rho _q})$ which verify the conditions $%
(a), $ $(b)$ and $(c)$ from $1.$ and there is a unitary operator $U_q$ in $%
L_{B_q}(E_{\rho _q,}F_q)$ which verifies the conditions $i),$ $ii)$ and $%
iii) $ from $2.$

Let $(\Phi _{\rho _{}};V_{\rho _{}};E_{\rho _{}})$ be the KSGNS\
construction associated with $\rho $. According to the proof of Theorem 4.6
in [\textbf{2}] , $(\pi _q)_{*}\circ \Phi _{\rho _{}}=\Phi _{\rho _q};$ $%
(\pi _q)_{*}(V_{\rho _{}})=V_{\rho _q};$ $(E_{\rho _{}})_q=E_{\rho _q}$ for
all $q\in S(B)$ and $\left( U_q\right) _q$ is a coherent sequence in $%
L_{B_q}(E_{\rho _q},F_q)$. It is not difficult to check that for each $g\in
G,$ $(v_g^{\rho _q})_q$ is a coherent sequence in $L_{B_q}(E_{\rho _q})$,
and the map $g\mapsto v_g^{\rho _{}},$ where $v_g^{\rho _{}}$ is an element
in $L_B(E_{\rho _{}})$ such that $\left( \pi _q\right) _{*}\left( v_g^{\rho
_{}}\right) =v_g^{\rho _q}$ for all $q\in S(B)$ is a unitary representation
of $G$ on $E_{\rho _{}}.$ Also it is not difficult to check that $(\Phi
_{\rho _{}},v^{\rho _{}},E_{\rho _{}})$ is a covariant representation of $%
(G,A,\alpha )$ which verifies the conditions $(a),$ $(b)$ and $(c)$ from $1.$

Let $U\in L_B(E_{\rho _{}},F)$ such that $\left( \pi _q\right) _{*}\left(
U\right) =U_q$ for all $q\in S(B).$ Clearly $U$ is a unitary operator in $%
L_B(E_{\rho _{}},F)$ and it verifies the conditions $i),$ $ii)$ and $iii)$
from $2.$ q.e.d.

\smallskip

\begin{Remark}
The covariant representation $(\Phi _{\rho _{{}}},v^{\rho _{{}}},E_{\rho
_{{}}})$ of $\left( G,A,\alpha \right) $ induced by $\rho $ is unique up to
unitary equivalence.
\end{Remark}

\smallskip

\smallskip From Proposition 3.4 and Theorem 3.6 we obtain the following
corollary.

\begin{Corollary}
Let $(G,A,\alpha )$ be a locally $C^{*}$-dynamical system such that $\alpha $
is an inverse limit action, let $u$ be a unitary representation of $G$ on a
Hilbert module $E$ over a locally $C^{*}$-algebra $B$, and let $\rho $ be a $%
u$ -covariant, non-degenerate, continuous completely positive linear map
from $A$ to $L_{B}(E)$. Then $\rho $ induces a non-degenerate representation
of the crossed product $A\times _{\alpha _{{}}}G$ on a Hilbert $B$ -module.
\end{Corollary}

\smallskip

The following proposition is a generalization of Proposition 2 in [\textbf{4}%
] in the context of locally $C^{*}$-algebras.

\begin{Proposition}
Let $(G,A,\alpha )$ be a locally $C^{*}$-dynamical system such that $\alpha $
is an inverse limit action, let $B$ be a locally $C^{*}$-algebra, let $E$ be
a Hilbert \ $B$ -module and let $u$ be a unitary representation of $G$ on \ $%
E.$ If $\rho $ is a $u$ -covariant, non-degenerate, continuous completely
positive linear map from $A$ to $L_{B}(E),$ then there is a unique
completely positive linear map $\varphi $ from $A\times _{\alpha ^{{}}}G$ to 
$L_{B}(E)$ such that 
\[
\varphi (f)=\int\limits_{G}\rho (f(g))u_{g}dg 
\]
for all $f\in C_{c}(G,A).$ Moreover, $\varphi $ is non-degenerate.
\end{Proposition}

PROOF. By Theorem 3.6 there is a covariant representation $(\Phi _{\rho
_{}},v^{\rho _{}},E_{\rho _{}})$ of $(G,A,\alpha )$ and an element $V_{\rho
_{}}$ in $L_B(E,E_{\rho _{}})$ such that $\rho (a)=V_{\rho _{}}^{*}\Phi
_{\rho _{}}(a)V_{\rho _{}}$ and $v_g^{\rho _{}}V_{\rho _{}}=V_{\rho _{}}u_g$
for all $a\in A,$ and for all $g\in G.$

Let $\Phi _{\rho _{}}\times v^{\rho _{}}$ be the representation of $A\times
_{\alpha ^{}}G$ associated with $(\Phi _{\rho _{}},v^{\rho _{}},E_{\rho
_{}}).$ We define $\varphi :A\times _{\alpha ^{}}G$ to $L_B(E)\;$by 
\[
\varphi (x)=V_{\rho _{}}^{*}(\Phi _{\rho _{}}\times v^{\rho _{}})(x)V_{\rho
_{}}. 
\]
Clearly $\varphi $ is a continuous completely positive linear map from $%
A\times _{\alpha ^{}}G$ to $L_B(E).$ Let $\{e_{\lambda _{}}\}_{\lambda \in
\Lambda }$ be an approximate unit for $A\times _{\alpha ^{}}G$ and let $\xi
\in E.$ Since $\Phi _{\rho _{}}\times v^{\rho _{}}$ is non-degenerate, by
Proposition 4.2 in [\textbf{2}] 
\[
\lim\limits_{\lambda _{}}V_{\rho _{}}^{*}(\Phi _{\rho _{}}\times v^{\rho
_{}})(e_{\lambda _{}})V_{\rho _{}}\xi =V_{\rho _{}}^{*}V_{\rho _{}}\xi =\xi
. 
\]
This implies that the net $\{\rho \left( e_{\lambda _{}}\right) \}_{\lambda
\in \Lambda }$ converges strictly to the identity map on $E,$ and so $%
\varphi $ is non-degenerate.

For $f\in C_c(G,A)$ we have 
\begin{eqnarray*}
\varphi (f) &=&V_{\rho _{}}^{*}(\Phi _{\rho _{}}\times v^{\rho
_{}})(f)V_{\rho _{}}=\int\limits_GV_{\rho _{}}^{*}\Phi _{\rho
_{}}(f(g))v_g^{\rho _{}}V_{\rho _{}}dg \\
&=&\int\limits_GV_{\rho _{}}^{*}\Phi _{\rho _{}}(f(g))V_{\rho
_{}}u_gdg=\int\limits_G\rho (f(g))u_gdg
\end{eqnarray*}
and since $C_c(G,A)$ is dense in $A\times _{\alpha ^{}}G,$ $\varphi $ is
unique with this property. q.e.d.

\smallskip

\begin{Corollary}
Let $(G,A,\alpha )$ be a locally $C^{*}$-dynamical system, let $B$ be a
locally $C^{*}$-algebra, let $E$ be a Hilbert \ $B$ -module, let $u$ be a
unitary representation of $G$ on \ $E$, and let $\rho $ be a $u$ -covariant,
non-degenerate, continuous completely positive linear map from $A$ to $%
L_{B}(E)$. If $G$ is a compact group, then there is a unique completely
positive linear map $\varphi $ from $A\times _{\alpha ^{{}}}G$ to $L_{B}(E)$
such that 
\[
\varphi (f)=\int\limits_{G}\rho (f(g))u_{g}dg 
\]
for all $f\in C_{c}(G,A).$ Moreover, $\varphi $ is non-degenerate.
\end{Corollary}

PROOF. The corollary follows from Proposition 3.9, since $G$ is compact and
then $\alpha $ is an inverse limit action of $G$ on $A$ [\textbf{12}, Lemma
5.2].q.e.d.

Department of Mathematics, Faculty of Chemistry, University of Bucharest,
Bd. Regina Elisabeta nr.4-12, Bucharest, Romania\ \ 

\ \ \ \ \ \ \ \ mjoita@fmi.unibuc.ro


\begin{thebibliography}{99}
\bibitem{}  M. Fragoulopoulou, \textit{An introduction to the representation
theory of topological }$*$\textit{\ -algebras, }Schriftenreihe, Univ.
M\"{u}nster, \textbf{48}(1988), 1-81.

\bibitem{}  M. Joita, \textit{Strict completely positive maps between
locally }$C^{*}$\textit{-algebras and representations on Hilbert modules},
J. London Math. Soc. (2), \textbf{66}(2002), 421-432.

\bibitem{}  M. Joita, \textit{Crossed products of locally }$C^{*}$\textit{%
-algebras}, Rocky Mountain J. Math. (to appear).

\bibitem{}  A. Kaplan, \textit{Covariant completely positive maps and
liftings, }Rocky Mountain J. Math. \textbf{23}(1993), 939-946.

\bibitem{}  G. G. Kasparov, \textit{Hilbert }$C^{*}$\textit{-modules:
Theorem of Stinespring and Voiculescu}, J. Operator Theory \textbf{4}(1980),
133-150.

\bibitem{}  E. C. Lance, \textit{Hilbert }$C^{*}$\textit{-modules. A toolkit
for operator algebraists, }London Mathematical Society Lecture Note Series
210, Cambridge University Press, Cambridge 1995.

\bibitem{}  A. Mallios, \textit{Topological algebras: Selected Topics},
North Holland, Amsterdam, 1986.

\bibitem{}  W. L. Paschke, \textit{Inner product modules over }$B^{*}$%
\textit{-algebras}, Trans. Amer. Math. Soc. \textbf{182}(1973), 443-468.

\bibitem{}  V. Paulsen, \textit{A covariant version of Ext,} Michigan, Math.
J. 29(1982), 131-142.

\bibitem{}  G. K. Pedersen, $C^{*}$-\textit{algebras and their automorphism
groups,} Academic Press, London, New-York, San Francisco, 1979.

\bibitem{}  N. C. Phillips, \textit{Inverse limit of }$C^{*}$\textit{%
-algebras,} J. Operator Theory, \textbf{19}(1988), 159-195.

\bibitem{}  N. C. Phillips, \textit{Representable }$K$\textit{-theory for }$%
\sigma $\textit{\ -}$C^{*}$\textit{-algebras},\textit{\ }$K$-Theory, \textbf{%
3}(1989),5, 441-478.

\bibitem{}  W. Stinespring, \textit{Positive functions on }$C^{*}$\textit{%
-algebras}, Proc. Amer. Math. Soc., \textbf{6}(1955), 211-216.
\end{thebibliography}
\end{document}